\documentclass{amsart}
\usepackage{graphicx}
\usepackage{graphicx,esint}
\usepackage{amssymb,amscd,amsthm,amsxtra}
\usepackage{latexsym}
\usepackage{epsfig}

\numberwithin{equation}{section}

\newtheorem{thm}{Theorem}[section]
\newtheorem{cor}[thm]{Corollary}
\newtheorem{lem}[thm]{Lemma}
\newtheorem{prop}[thm]{Proposition}
\theoremstyle{definition}

\theoremstyle{remark}

\numberwithin{equation}{section}

\newcommand{\R} {\mathbb R}
\newcommand{\C} {\mathbb C}

\newcommand{\Z} {\mathbb Z}
\newcommand{\disp} {\displaystyle}

\DeclareMathOperator{\re}{Re}
\DeclareMathOperator{\im}{Im}

\def\a{\alpha}
\def\b{\beta}
\def\g{\gamma}
\def\G{\Gamma}
\def\d{\delta}
\def\D{\Delta}
\def\e{\epsilon}
\def\f{\varphi}
\def\F{\Phi}

\def\l{\lambda}

\def\n{\nu}

\def\p{\psi}

\def\r{\rho}
\def\s{\sigma}
\def\th{\theta}

\def\x{\xi}

\def\y{\eta}
\def \z {\zeta}

\def\de{\partial}

\begin{document}

\pagestyle{plain}
\title{
Singularities of the wave trace \\
near cluster points of the length spectrum}
\author{Yves Colin de Verdi\`ere}
\address{Institut de Fourier, U. de Grenoble I}
\email{\tt yves.colin-de-verdiere@ujf-grenoble.fr}
\author{Victor Guillemin}
\address{Massachusetts Institute of Technology}
\email{\tt vwg@math.mit.edu}
\author{David Jerison}
\address{Massachusetts Institute of Technology}
\email{\tt jerison@math.mit.edu}
\thanks{supported in part by NSF grant DMS-0244991}
\maketitle

\section{Introduction}
Consider the Laplace operator $\D$ on the disk $D$,
\[
\D = \de^2_{x_1} + \de^2_{x_2}; \quad D = \{x\in \R^2: x_1^2 + x_2^2 < 1\}
\]
and let $u_j$ be the (real-valued, normalized) Dirichlet
eigenfunctions
\[
\D u_j = -\l_ju_j, x\in D; \quad u_j(x) = 0, \ x\in \de D; 
\quad \int_D u_j(x)^2 dx = 1
\]
with eigenvalues $0 < \l_1 < \l_2 \le \l_3 \le \cdots$, indexed with multiplicity.
The wave trace is the sum 
\[
h(t) = \sum_{j=1}^\infty  e^{i\sqrt{\l_j} t},
\]
which converges in the sense of distributions.  The purpose
of this article is to announce the following theorem.
\begin{thm}\label{theorem} $h(t)$, an infinitely differentiable function 
on $2\pi < t < 8$, has a finite limit and is infinitely differentiable
at $t = 2\pi$ from the right.
\end{thm}
The significance of $2\pi$ is that 
it is a cluster point of the length spectrum from the left ($t<2\pi$),
as described in more detail below.   It is easy to verify that
there are no geodesics of length in between $2\pi$ and $8$;
it follows then from \cite{GM} that $h(t)$ is smooth in $2\pi < t < 8$. 
The content of the theorem is that $h$ is smooth from the right
up to the endpoint $2\pi$.   The same proof applies to every
cluster point $2\pi \ell$ of the length spectrum of geodesic flow on
the disk.

We recall now the relationship between $h(t)$ and the wave equation.
Consider the initial value problem for the wave equation,
\begin{equation*}
\begin{aligned}
(\de_t^2 - \D )u(t,x) &= 0, \ x\in D, \ t>0 \\
u(t,x) &= 0, \ x\in \de D, \ t>0 \\
u(0,x) &= f(x), \ x \in D \\
\de_t u(0,x) &= g(x), \ x \in D 
\end{aligned}
\end{equation*}
The solution is 
\[
u(t,x) = \int_D K_t^{(1)}(x,y) f(y)dy + \int_D K_t^{(2)}(x,y)g(y)dy
\]
where
\[
K_t^{(1)}(x,y) = \sum_{j=1}^\infty \cos(\sqrt{\l_j} t) u_j(x)u_j(y);
\quad
K_t^{(2)}(x,y) = \sum_{j=1}^\infty 
\frac{\sin(\sqrt{\l_j} t)}{\sqrt{\l_j}} u_j(x)u_j(y)
\]
and the trace of the operator with kernel $K_t^{(1)}$ is 
\[
\int_D K_t^{(1)}(x,x)dx  
= \sum_{j=1}^\infty \cos(\sqrt{\l_j} t)  = \re h(t)
\]
The same proof that shows that $h(t)$ is smooth as $t\to (2\pi)^+$
also shows that the trace
\[
\int_D K_t^{(2)}(x,x)dx  = \sum_{j=1}^\infty \sin(\sqrt{\l_j} t)/\sqrt{\l_j}
\]
is smooth as $t\to (2\pi)^+$.

The close connection between between the length spectrum  and 
the singularities of $h(t)$ was discovered by way of spectral 
and inverse spectral problems.  In \cite{CdV},  Colin de Verdi\`ere showed 
that on a generic compact Riemannian manifold without boundary, the length spectrum
is determined by the spectrum (list of eigenvalues of the Laplacian with 
multiplicity).
Duistermaat and Guillemin \cite{DG} and Chazarain \cite{C}
showed that the singular support of $h(t)$ is contained in
the length spectrum and that the two sets are equal generically.  
Indeed, in the generic case, the singularity of $h(t)$ for 
$t$ near $L$, the length of a geodesic, resembles
a negative power of $t-L$, or, more precisely, a conormal distribution.
In case the boundary is non-empty, Andersson and Melrose \cite{AM} introduced 
the notion of generalized geodesic length spectrum and proved the 
analogous inclusion for the singular support of $h$.  
Subsequent work on inverse spectral problems for domains in the plane
and, more generally, for manifolds with boundary can be found in 
\cite{GM,CdV,Z,HeZ}.

The periodic geodesics on the disk (with reflecting boundary)
have lengths
\[
L_{k,\ell} = k (2\sin (\pi\ell/k), \ \ell = \pm 1, \, \pm 2, \dots, \ k = 2, \, 3, \, 4, \dots
\]
with $k$ the number of segments (or reflections) of the
trajectory and $\ell$ the winding number of the trajectory
around the origin. The degenerate cases, $(2,1)$, $(4,2)$, $(6,3)$, \dots,  
correspond
to the trajectory that traverses a diameter $2$, $4$, $6$, \dots times.
It follows from \cite{GM} that $h(t)$ is singular at 
$t=\pm L_{k,\ell}$ and smooth in the complement of the closure of these
points.  In particular, $h(t)$ is singular when $t$ is the circumference of each
regular polygon,
\[
L_{2,1} < L_{3,1} < L_{4,1} < \cdots < L_{k,1} \to 2\pi \quad k\to \infty
\]
and h(t) is smooth in $2\pi < t < 8$, since the shortest
periodic geodesic with length greater than $2\pi$ has length $8=L_{4,2}$ 
(the $2$-gon traced twice).

Thus the issue addressed here that is not addressed in previous works,
is the behavior of $h(t)$ near a cluster point of the length spectrum.
Although we examine only the case of the disk, which
is far from generic, we expect the analogue of Theorem
\ref{theorem} to be valid for any convex domain in place of the disk.  
See the final remarks, below.

It is reasonable
to conjecture that $h(t)$ has some power law (classical conormal) behavior 
as $t \to (2\pi)^+$.  What is surprising is that the power
is $0$ and $h(t)$ is infinitely differentiable.  The spikes at $t = L_{k,1}$ 
are asymmetrical and decay rapidly for $t> L_{k,1}$, so rapidly that
an infinite sum of them with singularities at points closer and closer
to $2\pi$ still converges at $t= 2\pi$ .   This is
explained on a technical level by the fact that $h(t)$ is represented
by sum of oscillatory integrals in which the first
derivative of the phase function tends to infinity along one ray.
Because the derivative of the phase is large, the phase 
changes quickly, and the corresponding integrals have more cancellation 
than one would obtain from a classical phase function.

\section{Outline of the Proof}

Our theorem is proved from systematic,  optimal symbol properties of the zeros of 
Bessel functions  or, equivalently, the eigenvalues of the Laplace operator for the
Dirichlet problem on the disk.  

Let $\r(m,n)$ denote the $m$th zero of the $n$th Bessel function $J_n$,
$\r(1,n) < \r(2,n) < \cdots$.  
In polar coordinates $x_1 = r\cos\th$, $x_2 = r\sin\th$, the eigenfunctions 
on the disk have the form 
\[
J_n(\r(m,n)r)\cos n\th , \quad J_n(\r(m,n)r)\sin n\th 
\]
with eigenvalue $-\l = -\r(m,n)^2$.  Since the multiplicity of
the eigenvalue is two when $n\ge 1$ and one when $n=0$, 
\[
h(t) = \sum_{(m,n)\in \Z^2} \p_1(m)\p_2(n) e^{it\r(m,n)}
\]
for any smooth cut-off functions $\p_1$ and $\p_2$ satisfying
\[
\p_1(m) = \begin{cases} 
1, \ m> 7/8 \\
0, \ m < 3/4
\end{cases}; 
\quad
\p_2(n) = 
\begin{cases}
0, \ n< -1/4 \\
1, \ -1/8<  n< 1/8 \\
2, \ 3/4 < n
\end{cases}
\]
Below we will extend $\r(m,n)$ in a natural way to be defined 
for real numbers $m$.  (The extension to real numbers $n$ will be the standard
one for Bessel functions.)  The Poisson summation formula then 
yields
\[
h(t) = \sum_{(k,\ell)\in \Z^2} h_{k,\ell}(t)
\]
where
\[
h_{k,\ell}(t) = \int_{\R^2} \p_1(m)\p_2(n) e^{it\r(m,n)- 2\pi i(km +\ell n)}dm\, dn
\]
Our main result, Theorem \ref{theorem}, follows immediately from
\begin{thm}\label{main est}
\[
\left| \frac{d^{N_1}}{dt^{N_1}} h_{k,\ell}(t)\right|
\le C_{N_1,N_2} (1 + |k| + |\ell|)^{-N_2}
\]
for $2\pi < t < 2\pi + 1/10$
\end{thm}

The domain of integration for $h_{k,\ell}$ is the quadrant 
$m \ge 3/4$, $n \ge -1/4$.  We will deduce the estimates on $h_{k,\ell}$ 
from symbol estimates for $\r(m,n)$.  

In the sector range of the parameters, $m\ge cn$, for any fixed $c>0$, 
the zeros of the Bessel functions satisfy ordinary symbol estimates
as follows.
\begin{prop}\label{rhoclassical}  Fix a constant $c_0>0$.  If $m\ge 3/4$, 
$n\ge -1/4$ and $m\ge c_0n$, then 

a) $\disp |\de_m^j \de_n^k \r(m,n)| \le C_{j,k} (m + n)^{1-j-k}$

b) $\disp
(\de_m \r,\de_n \r) = (\pi/\sin \a ,\a/\sin \a) + O((m + n)^{-1})$
as $(m,n) \to \infty$
where $\a $ is defined by\footnote{We extend the definition of $\a$ continuously
across $n=0$, $\a = \pi/2$ by the reciprocal equation $1/(\tan \a - \a) = n/\pi m$.}
\[
\tan \a - \a = \pi m/n
\]
\end{prop}
In the range of values of $(m,n)$ complementary to Proposition \ref{rhoclassical},
the appropriate symbol-type estimates involve fractional powers of $m$ and $n$,
and $\de_m\r$ does tend to infinity.  (For $m$ fixed $n\to \infty$
it turns out that $\de_m \r \approx n^{1/3}$.  For our purposes
the subtlest and most important bound will be the lower bound
on the size of $\de_n\r$.)
\begin{thm}\label{nonclass}  There is an absolute constant $c_1>0$ 
such that if $3/4 \le m \le  c_0 n$ ($c_0$  from Proposition \ref{rhoclassical})
then 

a)  $\disp |\de_m^j\de_n^k (\r(m,n)-n)| \le C_{j,k}m^{2/3 - j} (m + n)^{1/3 - k}$

b) $\disp \de_n(\r(m,n) -n) \ge c_1 m^{2/3} n^{-2/3}$
\end{thm}
Proposition \ref{rhoclassical} and Theorem \ref{nonclass} will
be proved using asymptotic expansions of Bessel functions,
which we derive by the method of Watson (\cite{W}, p. 251) 
starting from the Debye contour integral representation.
\begin{equation}\label{hint}
J_\n(x) + iY_\n(x) =
\frac{1}{\pi i} \int_{-\infty}^{\infty + \pi i} 
e^{x\sinh z - \n z} dz
\end{equation}
We cannot merely quote Watson's asymptotic
expansion because we need to differentiate it.  To some 
extent these differentiated estimates were carried
out already by Ionescu and Jerison \cite{IJ}, but we need quite
a bit more detailed asymptotics, especially in the transition region 
where $x-n \approx Cn^{1/3}$.  Also, one needs to choose
the right coordinate system since differentiation in some directions behaves
differently from others.

Theorem \ref{main est} is proved by integration by parts.  Consider
the phase function of $h_{k,\ell}(t)$,
\[
Q= i(t\r(m,n) - 2\pi km - 2\pi \ell n)
\]
Since $Q$ is smooth, the only issue is the asymptotic behavior
as $(m,n)$ tends to infinity.  It is easy to show
from Proposition \ref{rhoclassical} (b) that for
large $(m,n)$ the critical points of the phase $Q$ occur
near $t= L_{k,\ell}$ ($\a = \pi\ell/k$).  

The rest of the paper is organized as follows. In Section 3,
we carry out the proof of Theorem \ref{main est},
dividing the $(m,n)$ quadrant of integration into two sectors.  
In the sector where $m\to \infty$, Proposition \ref{rhoclassical} and
standard integration by parts and non-stationary phase methods apply.
In the sector where $n\to \infty$, $3/4 \le m \le c_0n$, 
the lower bound on $\de_n Q$ given by Theorem \ref{nonclass}
is used when integrating by parts with respect to $n$.  On the other hand,
when integrating by parts in the $m$ variable, we will use
the oscillation of $e^{ikm}$ only and include the rest of the factors
in the exponential $e^Q$ in the amplitude.

In Section 4, we prove the symbol estimates and asymptotic expansions
for Bessel function following the method of steepest descent (Debye contours).
In Section 5, we deduce the symbol estimates for the zeros of Bessel functions
stated above. We conclude with remarks about the relationship
with earlier work and about the methods that can be expected to lead
to the analogous result when the disk is replaced by a convex domain.

\section{Integration by Parts}

We now deduce our main estimate, Theorem \ref{main est}, from the
symbol estimates for $\r(m,n)$,  Proposition \ref{rhoclassical} and 
Theorem \ref{nonclass}.  
Denote
\[
Q= i(t\r(m,n) - 2\pi km - 2\pi \ell n)
\]
Then
\begin{equation}\label{intA}
e^Q = \frac{-i}{t\de_n \r - 2\pi \ell} \de_n e^Q 
\end{equation}
\begin{equation}\label{intB}
e^Q = e^{i(t\r - 2\pi \ell n} \frac{i}{2\pi k} \de_m e^{-2\pi i mk} 
\end{equation}

We will divide the region of integration into sectors.  Consider first
the ``non-classical'' region $3/4 \le m < c_0n$.  Define
\[
I_{k,\ell}(t) =
\int_{\R^2} \p_1(m)\p_1(cn/m)  e^{i(t\r(m,n) - 2\pi km - 2\pi \ell n)} dm \, dn
\]
We focus at first on the case $\ell = 1$, $Q= i(t\r(m,n) - 2\pi km - 2\pi n)$
from which all the singularities near $t= 2\pi$ arise.  
Applying formula \ref{intA} and integrating by parts,
\begin{align*}
I_{k,1}(t) 
&= \int \p_1(m)\p_1(cn/m)  \frac{-i}{t\de_n \r - 2\pi} \de_n e^Q dm \, dn  \\
& = 
\int \frac{-it \de_{nn} \r}{(t\de_n \r - 2\pi)^2}
\p_1(m)\p_2(cn/m) e^Q dm \, dn  \\
& \qquad + 
\int \frac{c}{m} 
\p_1(m)\p_1'(cn/m) \frac{i }{t\de_n \r - 2\pi}e^Q dm \, dn  
\end{align*}
Repeating we obtain 
\begin{align*}
I_{k,1}(t) 
&= c_1 
\int \frac{1}{m^2} \p_1(m)\p_1''(cn/m)  
\frac{1}{(t\de_n \r - 2\pi)^2} e^Q dm \, dn  \\
& \quad + 
c_2 \int \frac{1}{m} \p_1(m)\p_1'(cn/m)  
\frac{t\de_{nn}\r}{(t\de_n \r - 2\pi)^2} e^Qdm \, dn  \\
& \qquad + 
\int  \p_1(m)\p_1(cn/m)  
\left[\frac{c_3(\de^2_{n}\r)^2}{(t\de_n \r - 2\pi)^4}
+ \frac{c_4 t\de^3_{n}\r}{(t\de_n \r - 2\pi)^3}\right]
 e^Q dm \, dn  
\end{align*}
for appropriate coefficients $c_j(t)$ (polynomial in $t$).
After $N$ integrations by parts, $I_{k,1}$ is expressed as a linear
combination of terms with integrand
\[
e^Q 
\frac{(\de_n^2\r)^{a_2} (\de_n^3\r)^{a_3}
\cdots (\de_n^{N+1}\r)^{a_{N+1}}}
{(t\de_n \r - 2\pi)^{N+ a_2 + a_3 + \cdots + a_{N+1}}};
\qquad 
(a_2 + 2a_3 + \cdots + Na_{N+1} = N)
\]
times cutoff functions  $\p_1(m)\p_1(cn/m)$ or  derivatives of these
cutoff functions.   For $t\ge 2\pi$, the denominator has the
lower bound
\[
t\de_n \r - 2\pi \ge 2\pi(\de_n\r - 1) \ge c m^{2/3}n^{-2/3}
\]
from Theorem \ref{nonclass}(b).  Moreover, Theorem \ref{nonclass}(a)
says, in particular, that  
\[
|\de_n^j \r| \precsim m^{2/3}n^{1/3 -j}
\]
Denote $B = a_2 + a_3 + \cdots + a_{N+1}$. Then each of the 
integrands is bounded by
\begin{align*}
\left|\frac{(\de_n^2\r)^{a_2} (\de_n^3\r)^{a_3}
\cdots (\de_n^{N+1}\r)^{a_{N+1}}}
{(t\de_n \r - 2\pi)^{N+ a_2 + a_3 + \cdots + a_{N+1}}}\right| 
& \precsim 
\frac{(m^{\frac23}n^{\frac13-2} )^{a_2} (m^{\frac23}n^{\frac13 - 3})^{a_3}
\cdots (m^{\frac23}n^{\frac13 - (N+1)})^{a_{N+1}}}
{(m^{\frac23}n^{-\frac23})^{N+ a_2 + a_3 + \cdots + a_{N+1}}} \\
& = 
\frac{m^{\frac23 B}
n^{\frac13 B}
n^{-(2a_2 + 3a_3 + \cdots + (N+1)a_{N+1})}}
{(m^{2/3}n^{-2/3})^N
m^{\frac23 B}
n^{-\frac23 B}} \\
& 
= \frac{n^{-(a_2 + 2a_3 + \cdots + Na_{N+1})}}
{(m^{2/3}n^{-2/3})^N}
= \frac{n^{-N}}
{(m^{2/3}n^{-2/3})^N}   \\
& = m^{-\frac23 N}n^{-\frac13 N}
\end{align*}
In particular, $I_{k,1}(t)$ is represented by a convergent integral.

Next, to prove rapid decay in $k$ we integrate by parts in $m$ using
substitution \ref{intB}.  The first step is
\begin{align*}
I_{k,1}(t) 
&= 
\int e^{i(t\r - 2\pi n)}\frac{-i}{2\pi k} \de_m e^{-2\pi i mk} \p_1(m) \p_1(cn/m) 
dm \, dn \\
& = 
\int \frac{-t\de_m \r}{2\pi k}  e^{Q} \p_1(m) \p_1(cn/m) 
dm \, dn \\
& + 
\int \frac{i}{2\pi k} e^{Q} \p_1'(m) \p_1(cn/m) 
dm \, dn 
+ 
\int \frac{-cni}{2\pi m^2k} e^{Q} \p_1(m) \p_1'(cn/m) 
dm \, dn 
\end{align*}
The term with the factor $n/m^2$ also has $\p_1'(cn/m)$ so
that it is supported where $cn \approx m$ and the term $n/m^2 $
is comparable to $1/n \sim 1/m$.  All terms have a gain
of a factor $1/k$ except the one in which the derivative $\de_m$ falls 
on $\r$.  In that case, 
\[
|\de_m \r | \precsim n^{1/3} m^{-1/3}
\]
In all, one step of type (\ref{intB}) yields the factor
\[
\frac{n^{1/3}m^{-1/3}}{k}
\]
Now we consider systematically what happens when steps
of type (\ref{intB}) are applied after $N$ steps of type 
(\ref{intA}).  If the derivative $\d_m$ in the integration
by parts falls on a factor $\de_n^j \r$, then this
gets replaced by $\de_n^j \de_m\r/k$ and thus the bound is
improved by the very favorable factor 
\[
\frac{1}{mk}
\]
If the derivative $\de_m$ falls on the $e^{i(t\r - 2\pi n)}$ as in
the first step, we have, as before a factor
\[
\frac{n^{1/3}m^{-1/3}}{k}
\]
If the derivative falls on a cutoff, then the gain is $1/k$. 
Finally, if the derivative falls on the denominator $(t\de_n r - 2\pi)$
then it produces a factor
\[
|\frac{\de_n \de_m \r}{k(t\de_n\r - 2\pi)}| \precsim
\frac{m^{-1/3}n^{-2/3}}{km^{2/3}n^{-2/3}} = \frac{1}{km}
\]
In all, the worst case is the factor $n^{1/3}m^{-1/3}/k$ for
each integration by parts in $m$.  Thus if we integrate
by parts $N$ times in $n$ and $M$ times in $m$, the
integrand will be bounded by
\[
m^{-\frac23 N}n^{-\frac13 N} n^{\frac13 M}m^{-\frac13 M}k^{-M} 
\]
Therefore, if we choose $N$ sufficiently large that 
\[
\frac13 N - \frac13 M > 2
\]
then we obtain a convergent integrand that gives a bound on
the integral by $k^{-M}$.

For $\ell\neq1$, the bound is
much simpler.  $\de_n \r$ is very close to $1$ for small $m/n$, so
if $2\pi \le t \le 4\pi - \delta$ for any fixed $\delta >0$, then
the denominator in the integrands, 
\[
|t\de_n \r - 2\pi\ell| \approx |\ell|
\]
for all integers $\ell \neq 1$.  By integrating by parts $N$ times,
one finds the bound
\[
|I_{k,\ell}| \precsim \ell^{-N}
\]
for each $N$.  The rapid decrease in $k$ follows from similar reasoning
to that given above for $I_{k,1}$. 

Next, consider derivatives $(d/dt)^{N_1}I_{k,\ell}$.   The integral
representing this expression just has an extra factor of
$\r^{N_1}$ in the integrand.   This extra factor has size
 $(n^{1/3}m^{2/3})^{N_1}$ and symbol type bounds of the obvious kind
after differentiation with respect to $m$ and $n$.   Thus 
one can compensate for these higher powers by more integrations by
parts, and nearly the same proof as above shows that the derivatives 
of $I_{k,\ell}(t)$ are also rapidly decreasing in $(k,\ell)$.  This ends the main
portion of the proof.  

What remains is to make estimates for the integrand in
the region region $m > cn$ of integration. 
\begin{equation}\label{easysectorh}
h_{k,\ell}- 2I_{k,\ell}  =\int
[\p_1(m)\p_2(n)- 2\p_1(m)\p_1(cn/m)] e^{Q}dmdn
\end{equation}
For this region we use Proposition \ref{rhoclassical}.   
The informal idea is as follows.  Fix $\alpha>0$, and consider a ray in $(m,n)$
space defined by 
\[
\tan\a - \a = \pi m/n
\]
If $\nabla Q \to (0,0)$ as $(m,n)\to \infty$ along this ray, then
the asymptotic formula implies
\[
t\pi/\sin \a = 2\pi k; \quad t \a \sin \a = 2\pi \ell,
\]
Thus if $Q$ has a ``critical point near infinity''
we can solve these equations for $\a$ and $t$ and find
\[
\a = \pi \ell/k; \quad t = 2k \sin (2\pi \ell/k) = L_{k,\ell}
\]
This explains the singularities at $t= L_{k,\ell}$.  The
fact that the phase is nonstationary at other values of $t$
will lead to a proof that $h_{k,\ell}(t)$ is smooth at each point 
$t\neq L_{k,\ell}$.  

In more detail, first consider the range $|(m,n)| \le C$,
truncating the integrand with a smooth bump function in $(m,n)$ 
variables.  In that range, since the derivatives of $\rho$ are
bounded, 
\[
|\de_n Q|\ge c |\ell|, \quad |\de_m Q|\ge c |k|,
\]
for sufficiently large $|k|$ and $|\ell|$.  Hence, writing
\[
e^Q = (1/\de_n Q)\de_n e^Q, \quad
e^Q = (1/\de_m Q)\de_m e^Q,
\]
and integrating by parts, one finds that the integral
decays like $1/(|k| + |\ell|)$.  Repeating $N$ times, one 
finds that the integral decays like $O(1/(|k| + |\ell|)^N)$ for
any $N$. 

Next we turn to the range $|(m,n)| \ge C$, $m\ge cn$.  In this case,
we will also do integration by parts either in the variable $m$ or $n$,
and use lower bounds on $|\de_n Q|$ or $|\de_mQ|$.

Here we will use the asymptotic formula
for $(\de_m\r, \de_n \r)$ of Proposition \ref{rhoclassical} b
and restrict $t$ to $2\pi \le t \le 2\pi + \delta$ for suitable
small number $\delta$.   

First we confirm 
\begin{equation}\label{n-deriv}
|\de_n Q| = |t\a/\sin \a -2\pi \ell + O(1/(m+n))| \ge c (|\ell| +1),
\end{equation}
Recall that  $\a$ is defined 
by $\tan\a - \a = \pi m/n$.
If $n >0$, then, since $m \ge cn$, we have 
$\a_0 \le \a \le \pi/2$ for some fixed $\a_0>0 $ depending on $c$.  
In the remaining range, $0 \ge  n \ge -1/4$ and $m \ge C-1/4$,
which implies that $\pi/2 \le \a \le \pi/2 + \delta$
for some small $\delta$ of size on the order of $1/C$.  It
follows that 
\[
1 + \a_0^2/10 \le \a /\sin \a \le  \pi/2 + 4\delta
\]
The lower bound\footnote{The key here is that in the case $\ell =1$,  
$t\a\sin\a - 2\pi \ge 2\pi\a_0^2/10 >0$.  We have avoided the
critical points associated with $t = L_{k,1}\to 2\pi$.
They occur at infinity along a rays in $(m,n)$ space near the 
$n$ axis, rays that are not in the sector $m\ge cn$.}
on $\a/\sin\a$ bounds $|\de_nQ|$ from below
when $\ell\le 1$ and the upper bound on $\a/\sin \a$ 
bounds $|\de_nQ|$ from below in the case $\ell \ge 2$.  Thus
we have proved  \eqref{n-deriv}.  We use this bound and
integration by parts in $n$ in the range.
$|\ell| + C_2 \ge |k|$.

Lastly, if $k \ge |\ell| + C_2$, then (taking $C_2 = 10/\sin \a_0$)
we have 
\begin{equation}\label{m-deriv}
|\de_m Q| = |t\pi/\sin \a -2\pi k + O(1/(m+n))| \ge c (|k| + 1).
\end{equation}
On the other hand, if $k\le0$, then \eqref{m-deriv}
is obvious since $t\pi/\sin \a_0 \ge 2\pi^2$.  In the
range $|k| \ge |\ell| + C_2$, we use \eqref{m-deriv} and
integration by parts in the $m$ variable.  This concludes
the proof of Theorem \ref{main est}.

\section{Asymptotics of Bessel functions}

In this section, we establish the optimal symbol properties of 
\[
H_n(x):=J_n(x) + iY_n(x)
\]
as a functions of two variables $(x,n)$.  
(The function $H_\n(x)$, known as a 
Hankel function or Bessel function of the third kind, is denoted 
$H_\n^{(1)}(x)$ in Watson's treatise \cite{W} p. 73.)
The asymptotic formula for the Bessel functions as the order and
variable tend to infinity was discovered by Nicholson in 1910.  In
1918, Watson used the Debye contour representation to give an
appropriate bound on the error term.  In his treatise on Bessel functions
(\cite{W} p. 249), Watson says of his own method that
it is ``theoretically simple (though actually it is very laborious).''
To prove Theorem \ref{nonclass}, we will carry out the even more laborious
process of differentiating Watson's asymptotic formulas.

To state the symbol properties of $H_n(x)$ in 
the sector $n\le x \le 2n$, especially in the so-called transition 
region in which  $x$ is very close to $n$, will require a different 
coordinate system $(\b,\n)$.    For $x \ge n  > 0$, we write $\n = n$
and define $\b$ by $x\cos\b = n$.  Define $a(\n,\b)$ by 
\begin{equation}\label{adef}
H_\n(\n\sec\b) = e^{i\n(\tan\b - \b)} a(\n,\b)
\end{equation}
\begin{prop} \label{abound} 
If $\n \ge 1/2$ and $0 \le \b \le \pi/4$, and
$a(\n,\b)$ is defined by \eqref{adef}, then 
\[
|\de_\n^j \de_\b^k a| \precsim \n^{-1/2-j}\b^{-1/2 - k}
\]
\end{prop}
We will also need more detailed asymptotics involving
Airy functions.  Define the function $A(y)$ for $y\in \R$ as the 
solution to the equation 
\[
A''(y) + 2 y A(y) = 0
\]
with initial conditions
\[
A(0) = \G(1/3)6^{-2/3}(3 + i\sqrt3); \quad
A'(0) = \G(2/3) 6^{-1/3}(-3 + i\sqrt3)
\]
\begin{prop}\label{a-bbound}  Denote 
$\disp b(\n,\b) = e^{-i\n (\tan^3\b)/3} \n^{-1/3}A(y)$ with $\disp
y = (1/2)\n^{2/3}\tan^2\b$ and $a(\n,\b)$ from Proposition \ref{abound}.  Then for
$0 \le \b \le \pi/4$, 

a)  $|a-b| \precsim \n^{-1}$.

b) $|\de_\n(a-b)| \precsim \n^{-2}$.

c) $|\de_\b(a-b)| \precsim \n^{-2/3}$  [also  $\precsim 1/\n\b$ if $\b \ge \n^{-1/3}$]
\end{prop}

\begin{cor} \label{alowerbound} If $\n >> 1$ and $0 \le \b \le \pi/4$, and
$a(\n,\b)$ is defined by \eqref{adef}, then 

a) $\disp |a| \approx \n^{-1/2}\b^{-1/2}$ provided $\n^{-1/3} \le \b \le \pi/4$.

b) $\disp |a| \approx \n^{-1/3}$ provided $ \b \le \n^{-1/3}$.
\end{cor}

In the remainder of the section, we will prove Propositions \ref{abound} 
and \ref{a-bbound} and the corollary.  The range,
$x\ge 2n$ will be discussed at the end of the section.

\noindent  \emph {Proof of Proposition \ref{abound}}.  
Define the phase function $\f(z)$ by 
\[
\f(z) = \frac{1}{\n} (x\sinh z - \n z) = \sec \b \sinh z - z
\]
The contour integral \eqref{hint} can then be written
\[
H_\n(\n\sec\b) 
= \frac{1}{\pi i} \int_{-\infty}^{\infty + \pi i} e^{\n \f(z)} dz
\]
The contour of steepest descent\footnote{We follow \cite{W} p. 244 and 
pp. 249--252, except that where Watson uses $e^{-xr}$, we use $e^{-\n r}$ so
our expressions differ from his by factors $x/\n = \cos\b$.}
passes through $z = i\b$ and is
parametrized by the two curves, 
$z = \z_1(r,\b) + i\b$ and $z = \z_2(r,\b) + i\b$,
in which the functions
$\z_j(r,\b)$, $j=1,2$, solve
\begin{equation}\label{zetaeqn}
\f(\z + i\b) - \f(i\b) 
= (\sec\b)\sinh(\z + i\b) - \z - i\tan \b = -r , \quad r>0,
\end{equation}
and satisfy $\z_1(0,\b) = \z_2(0,\b) = 0$, 
$\re \z_1(r,\b) \le 0$, 
$\re \z_2(r,\b) \ge0$.   Moreover, $\z_1(r,\b) \to -\infty -i\b$ and
$\z_2(r,\b) \to \infty + i(\pi - \b)$ as $r\to\infty$.  Thus
\begin{align*}
H_\n(\n\sec\b) 
& = \frac{e^{\n \f(i\b)}}{\pi i} \int_{-\infty}^{\infty + \pi i} e^{\n (\f(z) - \f(i\b))} dz \\
&= \frac{e^{i\n (\tan\b - \b)}}{\pi i} 
\int_0^\infty e^{-\n r} 
(\de_r\z_2(r,\b)- \de_r \z_1(r,\b)) dr,
\end{align*}
and we have derived the formula for $a(\n,\b)$,
\begin{equation}\label{acontour}
a(\n,\b) = \frac{1}{\pi i} \int_0^\infty e^{-\n r} 
(\de_r\z_2(r,\b)- \de_r \z_1(r,\b)) dr
\end{equation}

\begin{lem}\label{zeta} For $j=1,2$, 
\[
|\de_r \de_\b^k \z_j(r,\b)| \precsim 
\begin{cases}  r^{-1/2}\b^{-1/2-k}, \quad & r < \b^3 \\
  r^{-(2+k)/3}, \quad & \b^3 < r   
\end{cases}
\]
\end{lem}
Lemma \ref{zeta} follows in a straightforward way 
from implicit differentiation and induction, but the proof takes
a few pages.  
Abbreviate by $\z$ the functions $\z_j(r,\b)$ 
along the contour of steepest descent satisfying (\ref{zetaeqn}).
Differentiating (\ref{zetaeqn}) with respect to $\b$ yields
\[
[\cosh(i\b + \z) - \cosh (i\b)]\de_\b \z 
+ (\tan\b) \sinh(\z + i\b) + i\cosh (\z+i\b) - i\sec \b = 0
\]
Therefore, using $\cosh(a+b) = \cosh a \cosh b + \sinh a \sinh b $
with $a = \z+i\b$ and $b= -i\b$, 
\[
\de_\b \z = \frac{-i(\cosh \z - 1)\sec\b}{\cosh(\z + i\b) - \cosh(i\b)}
\]
Define
\[
F_1(\z) = 1/\sinh(\z/2); \quad F_2(\z,\b) = 1/\sinh(\z/2 + i\b);
\quad G(\z,\b) = (\cosh \z -1)(\sec\b)/\z^2
\]
Since
\[
\cosh(i\b + \z) - \cosh (i\b) = 2\sinh(\z/2) \sinh(i\b + \z/2),
\]
we may write
\begin{equation}\label{zetabeta}
\de_\b \z = \frac{-i}{2} F_1 F_2 G \z^2
\end{equation}
Similarly,
\begin{equation}\label{zetar}
\de_r \z = \frac{-1}{2} F_1 F_2 \cos\b
\end{equation}
\begin{lem}\label{betainduction}  Along the contour of 
steepest descent $\z= \z_1(r,\b)$ and $\z = \z_2(r,\b)$, for all $k\ge0$,

a)  
$\disp |\de_\b^k F_1(\z(r,\b))| \precsim 
\begin{cases} |\z|^{-1}(\b + |\z|)^{-k} & \ \mbox{for} \quad |\z| \le 1\\
e^{-|\z|/2} & \ \mbox{for} \quad |\z| \ge 1
\end{cases} 
$

b)  
$\disp |\de_\b^k F_2(\z(r,\b),\b)| \precsim 
\begin{cases} (\b + |\z|)^{-k-1} & \ \mbox{for} \quad |\z| \le 1\\
e^{-|\z|/2} & \ \mbox{for} \quad |\z| \ge 1
\end{cases} 
$

c)  
$\disp |\de_\b^k G(\z(r,\b),\b)| 
\precsim 
\begin{cases} (\b + |\z|)^{-k} & \ \mbox{for} \quad |\z| \le 1\\
e^{|\z|}/|\z|^2 & \ \mbox{for} \quad |\z| \ge 1
\end{cases} 
$

d) 
$\disp |\de_\b^{k+1} \z(r,\b)| \precsim 
\begin{cases} |\z|(\b + |\z|)^{-k-1} & \ \mbox{for} \quad |\z| \le 1 \\
1 & \ \mbox{for} \quad |\z| \ge 1
\end{cases} 
$

\end{lem}
To begin the proof of Lemma \ref{betainduction}, denote 
by $D_\z F(\z,\b)$ the derivative of $F(\z,\b)$ with $\b$ fixed
and let $D_\b F(\z,\b) $ represent the derivative of $F(\z,\b)$ with $\z$ 
fixed.  (This is to distinguish from the partial derivative 
$\de_\b F(\z(r,\b),\b) $ representing the derivative with $r$ fixed.)
We bound $D_\z$ and $D_\b$ derivatives 
of $F_1$, $F_2$ and $G$ as follows.  On the curve of steepest descent 
($\z = \z_j(r,\b)$, $j=1,2$)
\begin{equation}\label{f1bound}
|D_\z^p F_1(\z)| = |(d/d\z)^pF_1(\z)|  \precsim 
\begin{cases} 
|\z|^{-k-1} & \quad \mbox{for} \ |\z|\le 1\\
e^{-|\z|/2} & \quad \mbox{for} \ |\z|\ge 1\\
\end{cases}
\end{equation}
\begin{equation}\label{f2bound}
|D_\z^p D_\b^sF_2(\z,\b)|
\precsim 
\begin{cases} 
(\b + |\z|)^{-p-s-1} & \quad \mbox{for} \ |\z|\le 1\\
e^{-|\z|/2} & \quad \mbox{for} \ |\z|\ge 1\\
\end{cases}
\end{equation}
\begin{equation}\label{gbound}
|D_\z^pD_\b^s G (\z,\b)|
\precsim 
\begin{cases} 
1 & \quad \mbox{for} \ |\z|\le 1\\
e^{|\z|}/|\z|^2 & \quad \mbox{for} \ |\z|\ge 1\\
\end{cases} 
\end{equation}
To prove (\ref{f1bound}) when $p=0$, note that if 
$\z = \x + i\y$ and $\x \ge 0$, then
\[
|\sinh(\z/2)| \ge (e^{\x/2}-1)/2 \ge 
\begin{cases} 
\x/4, & \quad  \x\ge0 \\
e^{\x/2}/4, & \quad  \x \ge 2 
\end{cases}
\]
Along the contour $\z = \x + i\y = \z_2(r,\b)$, $0 \le \y \le \pi - \b$.
(Recall that $0 \le \b \le \pi/2$.)  Furthermore, 
the slope 
of $\y$ as a function of $\x$ is at most $\sqrt3$ (\cite{W} 8.32, p.~240)
so that $\x \ge |\z| - \pi$ and $\x \ge |\z|/10$.  It follows that
$|F_1(\z)| = |1/\sinh(\z)| \le 40|\z|^{-1}$ for all $\z= \z_2(r,\b)$.
For $|\z| \ge 20$, we have $\x \ge 2$ and consequently
$\disp |F_1(z)| \le 4e^{-\x/2} \le 4e^{\pi/2} e^{-|\z|/2}$.  This proves
\eqref{f1bound} for $p=0$ and $\z = \z_2(r,\b)$.  The other branch $\z = \z_1$ is
similar, with the only difference that the curve $\z_1$ is in a horizontal
strip of the complex plane is below the $\x$-axis: $-\b \le \y = \im \z_1 \le 0$.
The case of $p=1,2,\dots$ are easy consequences of the same estimates.
A very similar proof to the one for (\ref{f1bound}) gives (\ref{f2bound}) 
because 
\[
|i\b + \z| \approx \b + |\z|
\]
along the contour.  The estimate (\ref{gbound}) is easy.

We can now carry out the proof of 
Lemma \ref{betainduction} by induction.  Parts (a)--(c)
for $k=0$ are the same as (\ref{f1bound}--\ref{gbound}) for
$p=0$.  Part (d) for $k=0$ ($k+1=1$) follows from (a)--(c) for $k=0$ 
and the formula (\ref{zetabeta}).

For the induction step, assume (a)--(d) are
valid up to $k$.  Part (a), $\de_\b^{k+1}F_1(\z(r,\b))$, is a linear 
combination of terms 
\[
(D_\z^pF_1)(\de_\b^{q_1}\z)(\de_\b^{q_2}\z) \cdots (\de_\b^{q_p}\z)
\quad (q_1 + \cdots + q_p = k+1; \ q_j\ge 1)
\]
Part (b), $\de_\b^{k+1} F_2(\z(r,\b),\b) $,
is a linear combination of terms
\[
(D_\z^p D_\b^s F_2(z,\b))
(\de_\b^{q_1}\z)(\de_\b^{q_2}\z) \cdots (\de_\b^{q_p}\z)
\quad (q_1 + \cdots + q_p + s = k+1; \ q_j\ge 1)
\]
Part (c),
$\de_\b^{k+1} G(\z(r,\b),\b) $,
is a linear combination of terms
\[
(D_\z^pD_\b^s G(z,\b))
(\de_\b^{q_1}\z)(\de_\b^{q_2}\z) \cdots (\de_\b^{q_p}\z)
\quad (q_1 + \cdots + q_p + s = k+1; \ q_j\ge 1)
\]
Since the induction hypothesis says that $\de_\b^q \z$ has the given bounds 
for all $1\le q\le k+1$,  each of the factors has appropriate bounds and 
multiplying them yields the appropriate bounds for (a)--(c) with $k$ 
replaced by $k+1$.  

Lastly, to prove the induction step for (d), observe that
\[
\de_\b^{k+1} (\de_\b \z) = \de_\b^{k+1}[(-i/2)F_1 F_2 G \z^2]
\]
is a linear combination of terms of the form
\[
(\de_\b^{k_1}F_1(\z))( \de_\b^{k_2}F_2(\z,\b))( \de_\b^{k_3}G)( \de_\b^{k_4}\z)( \de_\b^{k_5}\z )
\]
with $k_1 + k_2 + \cdots +k_5 = k+1$.  Since $k_j\le k+1$, we have already 
proved the appropriate bounds on each factor.  Multiplying them yields
the correct bound for $\de_\b^{k+2}\z$.  This ends the proof of 
Lemma \ref{betainduction}.

To convert the implicit bounds of Lemma \ref{betainduction} to ones in terms 
of $r$ and $\b$, observe that for $\z= \z_j(r,\b)$, 
\[
|\z| \approx \begin{cases} 
r^{1/2} \b^{-1/2}, & \quad r \le \b^3 \\
r^{1/3}, & \quad  \b^3\le r \le 1 \\
1 + \log r  , & \quad   1 \le r 
\end{cases}
\]
In the range $r\ge1$, there is also a more precise
estimate along the contour, namely,
\[  
e^{|\z|} \approx e^{|\x|} \approx r \qquad (\z = \x + i\y)
\]
It follows that
\[
|i\b + \z|\approx \b + |\z|  \approx \begin{cases} 
\b, & \quad r \le \b^3 \\
r^{1/3}, & \quad  \b^3\le r \le 1 \\
1 + \log r  , & \quad   1 \le r 
\end{cases}
\]

With these upper and lower bounds on $|\z|$ and $\b + |\z|$, one can rewrite
Lemma \ref{betainduction} as 
\begin{lem}\label{betainduction1}  Along the contour of 
steepest descent $\z= \z_1(r,\b)$ and $\z = \z_2(r,\b)$, for all $k\ge0$,

a)  
$\disp |\de_\b^k F_1(\z(r,\b))| \precsim 
\begin{cases} 
r^{-1/2}\b^{1/2 - k} & \quad r \le \b^3 \\
r^{-\frac13(k+1)} & \quad  \b^3 \le r \le 1 \\
r^{-1/2} & \quad  1 \le r
\end{cases} 
$

b)  
$\disp |\de_\b^k F_2(\z(r,\b),\b)| \precsim 
\begin{cases} 
\b^{-(k+1)} & \quad r \le \b^3 \\
r^{-\frac13(k+1)} & \quad \b^3 \le r \le 1 \\
r^{-1/2} & \quad \quad  1 \le r
\end{cases} 
$

c)  
$\disp |\de_\b^k G(\z(r,\b),\b)| 
\precsim 
\begin{cases} 
\b^{-k} & \quad  r \le \b^3 \\
r^{-k/3} & \quad \b^3 \le r \le 1 \\
r/(1 + \log r)^2 & \quad  \quad  1 \le r
\end{cases} 
$

d) 
$\disp |\de_\b^{k+1} \z(r,\b)| \precsim 
\begin{cases} 
r^{1/2}\b^{-3/2 -k} & \quad r \le \b^3 \\
r^{-\frac13 k} & \quad  \b^3 \le r \le 1 \\
1  & \quad  \quad  1 \le r
\end{cases} 
$
\end{lem}
It is now routine to confirm Lemma \ref{zeta}.  
Differentiate (\ref{zetar}),
\[
\de_\b^k\de_r \z = (-1/2)\de_\b^k (F_1 F_2 \cos\b)
\]
The right-hand side is a linear combination of terms of the form
\[
(\de_\b^{k_1}F_1 )(\de_\b^{k_2}F_2 )(\de_\b^{k_3}\cos \b),
\quad k_1 + k_2 + k_3 = k
\]
which are bounded using Lemma \ref{betainduction1}.  Now 
that Lemma \ref{zeta} is proved, Proposition
\ref{abound} follows from the formula 
\eqref{acontour} for $a(\n,\b)$.  

We now turn to the proof of Proposition \ref{a-bbound}.  
The function $\f$ from \eqref{zetaeqn} can be rewritten
\begin{equation}\label{zetaphase}
\f(z + i\b) - \f(i\b) = (i\tan \b)(\cosh z - 1)+  (\sinh z - z) 
\end{equation}
The Taylor approximations $\cosh z -1 \approx z^2/2$ and $\sinh z -z = z^3/6$
so the cubic approximation to \eqref{zetaphase} is 
\begin{equation}\label{Zphase} 
\F(z) = (i\tan \b)z^2/2 + z^3/6 
\end{equation}
Following Watson, we prove Proposition \ref{a-bbound} by 
comparing $\z_1(r,\b)$ and $\z_2(r,\b)$, the two solutions to 
\eqref{zetaeqn}, to the functions $Z_1(r,\b)$ and $Z_2(r,\b)$ solving 
the corresponding equation
\begin{equation}\label{Zeqn}
\F(Z) = -r, \quad r\ge0
\end{equation} 
with $Z_1(0,\b)= Z_2(0,\b) =0$, and $\re Z_1(r,\b) \le 0$, $\re Z_2(r,\b) \ge0$.

Following Watson again, curve $Z_1(r) \to -\infty -i\tan\b $ 
as $r\to\infty$, whereas
$Z_2(r)$ is asymptotic to the ray whose argument is $\pi/3$ as $r\to \infty$.
Thus we define
\[
b(\b,\n) = \int_{-\infty - i\tan\b}^{e^{i\pi /3}\infty} e^{\n \F(Z)} dZ
\]
In order to show that this function $b(\b,\n)$ is the same as the
function $b$ of Proposition \ref{a-bbound}, we evaluate
it using the contour consisting of the two rays
$Z = -i\tan \b - \xi$ and $Z = -i\tan\b + \xi e^{i\pi/3} $, $\xi \ge 0$,
rather than the steepest descent contour defined using $Z_1$ and $Z_2$.
With the parametrization given for these rays,
and the change of variable $\xi = u/\n^{1/3}$, one obtains
\begin{align}\label{Bformula}
b(\b,\n) &= 
e^{-i(\n \tan^3\b)/3}\int_0^\infty e^{-\n \xi^3/6} 
\left[\g e^{\g( \n \xi \tan^2\b)/2} + e^{-(\n \xi \tan^2\b)/2} \right]d\xi \\
&= e^{-i(\n \tan^3\b)/3}\n^{-1/3}A((\nu^{2/3}\tan^2 \b)/2) \qquad (\g = e^{\pi i/3}) \notag
\end{align}
where $A(t)$ is defined by 
\[
A(t) = \int_0^\infty e^{-u^3/6}[\g e^{\g t u} + e^{-tu}]du
\]
We claim that
\[
A''(t) + 2 t A(t) = 0
\]
Indeed, for $\l\in \C$, define
\[
F_\l(t) = \int_0^\infty e^{- u^3/6 + \l tu} du
\]
Then 
\begin{align*}
-1 &= \int_0^\infty (d/du) e^{- u^3/6 + \l tu} du  \\
&= \int_0^\infty (-u^2/2 + \l t) e^{-  u^3/6 + \l tu} du  \\
&= (-1/2\l^2)F_\l''(t)+ \l t F_\l(t)
\end{align*}
Hence,
\[
-\frac12 F_\l''(t) + \l^3 t F_\l(t)  = -\l^2
\]
Then $A(t) = \g F_{\g}(t) + F_{-1}(t)  $ and $\g^3 = -1$ give
the equation $A''(t) + 2 t A(t) = 0$, as desired.

Thus $A(t)$ is an Airy-type function, identified uniquely by
its value and derivative
\[
A(0) = \G(1/3)6^{-2/3}(3 + i\sqrt3); \quad
A'(0) = \G(2/3) 6^{-1/3}(-3 + i\sqrt3)
\]
Writing $A$ in terms of its real and imaginary parts, $A(t) = u(t) + iv(t)$,
we find that the Wronskian takes the constant value
\begin{equation} \label{wronskian}
u(t) v'(t) - u'(t) v(t) = 
u(0) v'(0) - u'(0) v(0) = \G(1/3)\G(2/3)\sqrt3 = 2\pi
\end{equation}

Having identified $b$ with the function of Proposition \ref{a-bbound},
we can now proceed with the proof.

\begin{lem}\label{Zest} For $j = 1,\ 2$, $0 \le \b \le \pi/4$,
\[
|\z_j| + |Z_j| \precsim
\begin{cases} 
r^{1/2}\b^{-1/2} &\quad r < \b^3 \\
r^{1/3} &\quad  r > \b^3 
\end{cases}
\]
\end{lem}
\noindent
Lemma \ref{Zest} is routine and the proof is omitted. 
In fact $\z_j$ grows more slowly than $Z_j$ for 
for $r>>1$ (like $\log r$), but we don't 
make use of this.

\begin{lem}\label{denom} 
For $j = 1,\ 2$, 

a) $\disp 
|\f'(\z_j)| = |i\tan \b \sinh \z_j + \cosh \z_j - 1| 
\succsim
\begin{cases} 
r^{1/2}\b^{1/2} &\quad r < \b^3 \\
r^{2/3} &\quad  \b^3 \le r \le  1 \\
r &\quad   1 \le r\\
\end{cases}
$

b) $\disp |\F'(Z_j)| = |i\tan \b Z_j + Z_j^2/2| 
\succsim
\begin{cases} 
r^{1/2}\b^{1/2} &\quad r < \b^3 \\
r^{2/3} &\quad  r > \b^3 
\end{cases}$

\end{lem}
Part (a) of this lemma is proved in  \cite{IJ} 9.15 and 9.16, p.~1072.  
The proof of part (b) is similar and is omitted.

\begin{lem}\label{zminusZ} For $j = 1,\ 2$, $0 \le \b \le \pi/4$,
\[
|\z_j -Z_j| 
\le 
\begin{cases} 
r^{3/2}\b^{-3/2} &\quad r \le \b^3 \\
r &\quad   \b^3 \le r \le 1 \\
r^{1/3} &\quad   1 \le r\\
\end{cases} 
\]
\end{lem}
\noindent \emph{Proof}.   Fix a small number $r_0>0$.  For $0r \ge r_0$, 
the estimate follows immediately from Lemma \ref{Zest}.  For this
argument we will use subscripts on the constants, because
some will depend on others.   For $|z|\le 1$,  
$|\f(z) - \F(z)| \le C_1(\b |z|^4 +|z|^5 )$.  
Fix $r$, $\b$ and $j=1$ or $2$.   We consider first
the case $0< r \le \b^3$. 
Denote $w_j = \z_j(r,\b) - Z_j(r,\b)$ and $Z_j= Z_j(r,\b)$.
Our goal is to show that 
\[
|w_j| \le Cr^{3/2}\b^{-3/2}
\]
The root in the appropriate half-plane of 
\[
\f(w + Z_j) + r = 0
\]
is $w= w_j$.  We will show that for a suitable constant $C$
to be chosen later, the curve 
\[
S = \{\f(w+Z_j)+r: |w| = Cr^{3/2}\b^{-3/2}\}
\]
encloses the origin.  Thus the root $w_j$ is inside.   We
will show that $S$ enclosed the origin by showing that
$\f(w+Z_j) + r - \F'(Z_j)w$ is suitably small.  
Recall
that $|\F'(Z_j)|\ge c_1 r^{1/2}\b^{1/2}$.  It follows that the circle
\[
S_0= \{\F'(Z_j)w: |w| = Cr^{3/2}\b^{-3/2}\}
\]
has radius at least $Cc_1r^2/\b$.

By Lemma \ref{Zest}, $|Z_j| \le C_2r^{1/2}\b^{-1/2} $.  We will
require that 
\begin{equation}\label{req1}
Cr^{3/2}\b^{-3/2} \le C_2 r^{1/2}\b^{-1/2} 
\end{equation}
so that $|w + Z_j| \le 2C_2r^{1/2}\b^{-1/2}$ and 
\[
|\f(w+Z_j)  - \F(w+Z_j)| \le C_1(\b|w + Z_j|^4 + |w + Z_j|^5) 
\le C_3  r^{2}/\b
\]
with a constant $C_3$ depending only on $C_1$ and $C_2$ (using
$r \le \b^3$).  Furthermore, since $F(Z_j) = -r$, 
\[
\F(w + Z_j) = -r + \F'(Z_j)w + (1/2)(i\tan \b + Z_j)w^2 + w^3/6
\]
and 
\[
|(i\tan \b + Z_j)w^2/2 + w^3/6| \le  
C^2 \b^{-2}r^3 + 
C_2 C^2 r^{7/2} \b^{-7/2}  + C^3r^{9/2} \b^{-9/2} 
\]
We will require
\begin{equation}\label{req2}
C^2 \b^{-2}r^3 + 
C_2 C^2 r^{7/2} \b^{-7/2}  + C^3r^{9/2} \b^{-9/2} \le (1/100)c_1Cr^2/\b
\end{equation}
Now our final requirement on $C$ is the lower bound
\begin{equation}\label{req3}
C_3 \le (1/4) c_1 C
\end{equation}
Combining these estimates, we have for $|w| = Cr^{3/2}\b^{-3/2}$, 
\[
|\f(w + Z_j) + r - \F'(Z_j)w| \le C_3 r^2/\b + (1/100)c_1Cr^2/\b
\le (1/2) c_1Cr^2/\b
\]
so that
\[
|\f(w + Z_j) + r| \ge (1/2) c_1Cr^2/\b.
\]
On the other hand $|\F'(Z_j)w| \ge c_1Cr^2/\b$.  So $S$ is a
loop surrounding the origin. 

Finally, we check that
all three requirements on $C$ are satisfied.   For \eqref{req3} fix
$C = 4C_3/c_1$.  Now that $C$ is fixed, we choose $r_0>0$ sufficiently
small that the other two inequalities \eqref{req1} and \eqref{req2} are 
satisfied for all $r$, $0< r\le r_0$, $0 < r \le \b^3$ ($0 \le \b \le \pi/4$).
The remaining case, $\b^3 \le r \le r_0$ is similar and slightly
simpler.  It will be omitted.   This concludes  the proof of Lemma \ref{zminusZ}.

\begin{lem}\label{rderiv} For $j = 1,\ 2$, 

a) $\disp |\de_r Z_j|  \le 
\begin{cases} 
r^{-1/2}\b^{-1/2} &\quad r < \b^3 \\
r^{-2/3} &\quad  r > \b^3 
\end{cases}
$

b) $\disp |\de_r\z_j| \le 
\begin{cases} 
r^{-1/2}\b^{-1/2} &\quad r < \b^3 \\
r^{-2/3} &\quad  \b^3 \le r \le  1 \\
r^{-1} &\quad   1 \le r\\
\end{cases}$

c) $\disp 
|\de_r(\z_j - Z_j)| 
\le 1 \quad \mbox{all} \ r
$
\end{lem}
\noindent \emph{Proof}.  Differentiating \eqref{Zeqn}, we find
\[
\F'(Z_j) \de_r Z_j = -1; \quad \f'(\z_j) \de_r \z_j = -1
\]
and the bounds of Lemma \ref{denom} imply  parts (a) and (b).  For part (c),
\begin{align*}
\de_r(Z_j-\z_j) 
&= \frac{\F'(Z_j) - \f'(\z_j)}{\F'(Z_j)  \f'(\z_j)} \\
&= \frac{\F'(\z_j) - \f'(\z_j)}{\F'(Z_j)  \f'(\z_j)}
+ \frac{\F'(Z_j) - \F'(\z_j)}{\F'(Z_j)  \f'(\z_j)} 
\end{align*}
Lemma \ref{denom} implies 
\[
|{\F'(Z_j)  \f'(\z_j)} | \succsim \max(r\b,r^{4/3})
\]
The formula of $\F'$ and Lemmas \ref{Zest} and \ref{zminusZ} imply
\[
|\F'(Z_j) - \F'(\z_j)| 
\precsim  (\b + |Z_j|+|\z_j|)|Z_j-\z_j| 
\precsim  
\begin{cases} 
r^{3/2}\b^{-1/2} &\quad r < \b^3 \\
r^{4/3} &\quad  \b^3 \le r \le  1 \\
r^{2/3} &\quad   1 \le r\\
\end{cases}
\]
Hence, $|(\F'(Z_j) - \F'(\z_j))/\F'(Z_j)  \f'(\z_j))| \precsim 1$
Lemma \ref{Zest} implies
\[
|\F'(\z_j) - \f'(\z_j)| \precsim \b |\z_j|^3 + |\z_j|^4
\precsim
\begin{cases} 
r^{3/2}\b^{-1/2} &\quad r < \b^3 \\
r^{4/3} &\quad  \b^3 \le r  
\end{cases}
\]
Hence, similarly, $|(\F'(\z_j) - \f'(\z_j))/\F'(Z_j)  \f'(\z_j)| \precsim 1$.
This concludes Lemma \ref{rderiv}.

\begin{lem}\label{rbeta} For $j = 1,\ 2$, 

a) $\disp |\de_r \de_\b Z_j|  \le 
\begin{cases} 
r^{-1/2}\b^{-3/2} &\quad r \le \b^3 \\
r^{-1} &\quad  \b^3 \le r
\end{cases}$

b) $\disp 
|\de_r\de_\b\z_j| \le 
\begin{cases} 
r^{-1/2}\b^{-3/2} &\quad r < \b^3 \\
r^{-1} &\quad  \b^3 \le r 
\end{cases} $

c)  $\disp |\de_r\de_\b(\z_j - Z_j)| \precsim r^{-1/3}$
\end{lem}
\noindent \emph{Proof}.   Differentiate the implicit equation 
with respect to $\b$ to obtain
\[
\de_\b Z_j = -(i\sec^2\b) Z_j^2/2\F'(Z_j)
\]
and 
\[
\de_\b \de_rZ_j = -(i\sec^2\b)\de_r Z_j[2Z_j \F'(Z_j) - F''(Z_j)Z_j^2]/2F'(Z_j)^2
\]
We have already bounded each of these terms and the bounds combine
to give Lemma \ref{rbeta} (a).  To prove (b), differentiate \eqref{zetaeqn} with respect to
$\b$ to obtain
\[
\de_\b \z_j = -i(\sec^2\b)(\cosh \z_j - 1)/\f'(\z_j)
\]
Then, differentiating with respect to $r$, 
\[
\de_r \de_b \z_j = -i(\sec^2\b) \de_r \z_j[ \sinh \z_j \f'(\z_j) - \f''(\z_j)(\cosh \z_j -1)]/\f'(\z_j)^2
\]
The estimates above for $\z_j$, $\de_r\z_j$, $\f'(\z_j)$, and 
\[
\f''(\z_j) \precsim
\begin{cases}
\b  & \quad  r \le \b^3 \\
r^{1/3}  &\quad  \b^3 \le r  \le 1\\
r  & \quad  1 \le r
\end{cases}
\]
combine to give part (b) of Lemma \ref{rbeta}.

To prove (c) write
\begin{align}\label{betadiff}
\de_\b(\z_j - Z_j) 
& =
-i\sec^2\b [Z_j^2(\f'(\z_j) - \f'(Z_j)) + Z_j^2(\f'(Z_j) - \F'(Z_j)) + \\
& \phantom{-i\sec^2\b [Z_j^2}
(Z_j - \z_j)^2\F'(Z_j) - 2(\cosh \z_j - 1 - \z_j^2/2) \F'(Z_j)]/\F'(Zj) \f'(\z_j) 
\end{align}
With $z$ is a point on the line segment from $Z_j$ to $\z_j$, 
\begin{align}\label{betadiffbound}
|\de_\b(\z_j - Z_j) |
& \le 
C|Z_j|^2|\f''(z)||Z_j-\z_j| + |Z_j|^2(\b|Z_j|^3 +  \\
 & \phantom{ccccc}
 |Z_j|^4) + |Z_j||Z_j-\z_j||\F'(Z_j)| + |\z_j|^4||\F'(Z_j)|]/
|\F'(Z_j)\f'(\z_j)| 
\end{align}
Using the preceding bounds and \eqref{betadiffbound},
\[
|\de_\b(\z_j - Z_j) |
\precsim
\begin{cases}
r^{3/2}\b^{-5/2} & \quad  r \le \b^3 \\
r^{2/3} &\quad  \b^3 \le r  \le 1
\end{cases}
\]
In the range $ 0 < r \le 1$, differentiation of \eqref{betadiff}
with respect to $r$ 
replaces terms $Z_j$ (with the bound $r^{1/2}\b^{-1/2}$) by $\de_r Z_j$ (with the
bound $r^{-1/2}\b^{-1/2}$) or the similar replacement of $\z_j$ with  $\de_r \z_j$.  
This results in a bound of the same type as \eqref{betadiffbound} with an extra factor of $1/r$.
In other words,
\[
|\de_r\de_\b(\z_j - Z_j)| 
\precsim
\begin{cases}
r^{1/2} \b^{-5/2} & \quad  r \le \b^3 \\
r^{-1/3} &\quad  \b^3 \le r  \le 1 \\
r^{-1} & \quad  1 \le r
\end{cases} 
\]
(The case $r\ge1$ follows separately from parts (a) and (b).)
In all three cases this is less than $r^{-1/3}$, so this  concludes Lemma \ref{rbeta}.

Define
\[
\quad
b(\b,\n) = \int_0^\infty e^{-\n r} 
(\de_r Z_2(r,\b)- \de_r Z_1(r,\b)) dr
\]
Then Lemmas \ref{rderiv} and \ref{rbeta} imply that 
\begin{align}
|a(\b,\n) - b(\b,\n)| &\le \int_0^\infty e^{-\n r} dr \le 1/\n \label{error1}\\
|\de_\n[a(\b,\n) - b(\b,\n)]| &\le \int_0^\infty r e^{-\n r} dr \le 1/\n^2
\label{error2} \\
|\de_\b[a(\b,\n) - b(\b,\n)]| &\le \int_0^\infty r^{-1/3} e^{-\n r} dr \le  \n^{-2/3} \quad
\label{error3}
\end{align}
This concludes the proof of Proposition \ref{a-bbound}.

Next, in order to deduce Corollary \ref{alowerbound}, we will prove some
estimates for $b(\n,\b)$ that will also be needed in the next section.
\begin{prop} \label{b-bound} If $y\ge 0$ and $\n \ge 1/2$, $0\le \b \le \pi/4$, then

a) $\disp |b| \approx 
\begin{cases} 
\n^{-1/3} & \   \b \le \n^{-1/3} \\
\n^{-1/2} \b^{-1/2} & \b \ge \n^{-1/3}
\end{cases}$

b) $\disp |\de_\n b| \precsim 
\begin{cases} 
\n^{-4/3} & \   \b \le \n^{-1/3} \\
\n^{-1/2} \b^{5/2} & \b \ge \n^{-1/3}
\end{cases}$

c) $\disp |\de_\b b| \precsim 
\begin{cases} 
\n^{1/3}\b & \   \b \le \n^{-1/3} \\
\n^{1/2} \b^{3/2} & \b \ge \n^{-1/3}
\end{cases}$
\end{prop}
\noindent 
{\em Proof}.  The proposition will 
follow easily from the formula for $b$ in terms of $A(y)$ and
the estimates for $y\ge 0$,
\begin{equation} \label{Airybound}
|A(y)| \approx (1 + y)^{-1/4}; \qquad |A'(y)| \approx (1 + y)^{1/4}
\end{equation}
\eqref{Airybound} is  well known, but we include a sketch of a proof.
The upper bounds are standard; indeed, the asymptotic behavior
as $y\to \infty$ follows from the fact that $A(y)$ is a multiple of the Hankel
function $\disp H_{1/3}((2y)^{3/2}/3)$ \cite{W} p. 252.  The lower bounds (for
all $y\ge0$) follow from the upper bounds and the fact that
the Wronskian \eqref{wronskian} is constant. 

Next, using \eqref{Airybound} we deduce (a).  Recall that  from \eqref{Bformula}, 
\[
|b(\n,\b)| = \n^{-1/3} |A(y)| \approx \n^{-1/3} (1+|y|)^{-1/4}; \quad y = (1/2)\n^{2/3}\tan^2\b
\]
When $\b \le \n^{-1/3}$, $|y| \precsim 1$ and $|b|\approx \n^{-1/3}$.  When $\b \ge \n^{-1/3}$,
\[
|b(\n,\b)| \approx \n^{-1/3} (\nu^{2/3}\b^2)^{-1/4}\approx \n^{-1/2} \b^{-1/2}
\]
Parts (b) and (c) of Proposition \ref{b-bound} are proved as follows.
Differentiating (\ref{Bformula}) gives
\[
|\de_\n b(\n,\b) | \precsim  \b^3\n^{-1/3}|A(t)|  + \n^{-2/3}\b^2 |A'(t)| 
\precsim \min(\n^{-1/2}\b^{5/2} , \n^{-4/3})
\]
and
\[
|\de_\b b(\n,\b) | \precsim  \n^{2/3}\b^2 |A(t)| + \n^{1/3}\b|A'(t)| 
\precsim  \min(\n^{1/2} \b^{3/2},  \n^{1/3} \b)
\]

Proposition \ref{b-bound} (a) and (\ref{error1})  imply
Corollary \ref{alowerbound}.

Finally, we discuss the range $\pi/4 \le \b < \pi/2$ and
beyond. To formulate this we return to the variables $(x,n)$.
Fix $c>0$ and let $\b_0>0$ be the smallest number
such that $\cos\b_0 = 1/(1+c)$.  Let $\b_1$, $\pi/2 < \b_1 <\pi$ satisfy
$\cos \b_1  = -1/4$.  For $x\ge1$, $n\ge -1/4$ and $x \ge (1+c)n$, 
define $\b(x,n)$ as the unique number $\b_0 \le \b  \le \b_1$ such 
that $x \cos \b = n$.  
\begin{prop}\label{atildebound} Let $x\ge 1$, $n\ge -1/4$ and $x \ge (1+c)n$
for a fixed $c>0$.  Define $\tilde a(x,n)$ by
\[
H_n(x) =  e^{i(x\sin\b - n\b)}\tilde a(x,n).
\]
with $\b= \b(x,n)$ defined above. Then
\[
|\de_x^j \de_n^k \tilde a(x,n)| \precsim x^{-1/2 - j - k}
\]
\end{prop}
To explain the connection with the previous notation, if
$\n = n$ and $x=  \n \sec \b$,  then $a(\n,\b) = \tilde a(x,n)$.
The distinction between the coordinate
systems is that $\de_\n$ is  the derivative with $\b$ 
(or equivalently $x/\n$ fixed), whereas $\de_n$ represents
the derivative with $x$ held fixed. The
$\de_\n$ direction is special when $\b$ is near $0$
and the $\de_n$ direction is special when $\b$ is near 
$\pi/2$.  The estimates we carried out in the range
u$0\le \b \le \pi/4$ can be extended to $\b \to \pi/2$,
but they require additional factors of $\sec\b$ which
tends to infinity.  They do not suffice: In the eventual analysis of
the behavior of zeros $\r(m,n)$, the range $x= \r(m,n) \ge (1+c)n$ 
corresponds to $m \ge cn$ and estimates in the $(\n,\b)$ coordinates
give rise to error terms of size $O(1/n)$ when what is needed is 
$O(1/m) = O(1/\r)$.

Proposition \ref{atildebound} was already proved in 
the case $n \ge 0$,  $j=0$,
$k=0,1,2$ in Theorem 9.1 (i) of \cite{IJ}.  
(In the notation $a_x(\n)$ of \cite{IJ}, $\n = n$,
$a_x(n)x^{-1/4}(x-n)^{-1/4} = \tilde a(x,n)$, and in
the range of variables specified here, $x-n\approx x$.)
The full proof of Proposition \ref{atildebound}
follows the same procedure as in \cite{IJ}
pp. 1068--1072, with the only extra ingredient
being the systematic treatment of derivatives of all orders,
which was already carried out above in the very similar proof 
of the symbol estimates for $a(\n,\b)$ in Proposition \ref{abound} (a).  
These details will be omitted.  We call attention to one difference. 
The integrals in the proof of Proposition \ref{abound} involve $e^{-\n r}dr$
as $\n\to \infty$, whereas in the proof of Proposition \ref{atildebound},
(following \cite{IJ}) the integrals involve $e^{-xr}dr$ and $x\to \infty$.  
We mention this in order to explain why the proof is unchanged when the 
range of $n$ is extended 
from $n=\n \ge 0$ to $n\ge -1/4$.  The range $n \le 0$ would
be problematic for integrals on $0< r < \infty$ involving $e^{-nr}dr$, 
but these integrands are replaced by ones involving $e^{-xr}dr$ with $x\ge 1$.

\section{Asymptotics of $\arg J_n+iY_n$ and of the zeros of Bessel functions}

Denote $H_n(x) = J_n(x) + iY_n(x)$. It is well known that for all real numbers $n$,
$|H_n(0^+)| = \infty$ and $|H_n(x)|$ is a decreasing function of $x$ for 
$x>0$ (Watson \cite{W} p. 446).  Moreover,
\[
\de_x \arg(H_n(x)) = \im \frac{H'_n(x)}{H_n(x)} = 
\frac{J_n'(x)Y_n(x) - Y_n'(x)J_n(x)}
{J_n(x)^2 + Y_n(x)^2} = \frac{2}{\pi x |H_n(x)|^2}>0
\]
Therefore, for all real values of $n$, as $x>0$ increases, 
$H_n(x)$ traces a simple spiral 
counterclockwise in the complex plane.  To choose a well-defined branch 
of the argument
\[
\th(x,n) = \arg(H_n(x))
\]
note that 
$Y_n(x) = (J_n(x)\cos(n\pi) - J_{-n}(x))/\sin(n\pi)$ and one has the asymptotic
formula $J_n(x) \sim (x/2)^n/\G(n+1)$ as $x\to 0^+$.  We deduce that for $n\ge0$,
\[
H_n(x)/|H_n(x)| \to -i = e^{-i\pi/2} \quad \mbox{as} \ x\to 0^+
\]
whereas for $0 \le n \le -1/2$, 
\[
H_n(x)/|H_n(x)| \to -\sin n\pi -i \cos n\pi = e^{-i(n+1/2) \pi} \quad \mbox{as} \ x\to 0^+
\]
Therefore, for $x>0$ sufficiently small $H_n(x)$ is in the 4th quadrant.
and a consistent definition of the branch is given by 
\[
\th(0^+,n) = -\pi/2 \quad (n\ge0)
\] 
and 
\[
\th(0^+,n) =  -(n+1/2)\pi \quad 0 \ge n \ge -1/2
\]
It then follows that the $m$th positive zero of $J_n(x)$ satisfies
\begin{equation}\label{rhodef}
\th(\r(m,n),n)= m\pi - \pi/2 
\end{equation}
for $m= 1, 2, \dots$.  This implicit equation can then be used to extend
the definition of $\r(m,n)$ for all real $m$ and $n$ satisfying $m+n>0$. 

Next, we establish a few preliminary upper and lower bounds for $\r(m,n)$.
It follows from the fact that $\de_x\th(x,n)>0$ that $\r(m,n)$ 
well-defined and infinitely differentiable.  In the range 
$|n|\le 1/4$,  the formula for $\de_x\th(x,n)$ and the 
estimates $|H_n(x)| \precsim  x^{-1/4}$ in $0 < x < 1$ and
$|H_n(x)| \precsim  x^{-1/2}$ in $1 \le x <\infty$ imply that 
\begin{equation}\label{easy1}
\r(m,n) \approx m \quad (m\ge 1/2, \quad |n|\le 1/4)
\end{equation}
In the range $n\ge 1/4$, $m \ge 3/4$,  we prove (well-known) upper
and lower bounds 
\begin{equation}\label{easy2}
\r(m,n) - n \approx m + m^{2/3}n^{1/3} 
\end{equation}
Let $n>0$, then according to \cite{W} (p. 485--487),
$n <  y_n$ where $y_n= \r(1/2,n)$ is the smallest positive 
zero of $Y_n(x)$.  Since $\th(y_n,n) = 0$, 
\[
m\pi - \pi/2 = \int_{y_n}^{\r(m,n)}  \de_x \th(x,n)\, dx
\]
For all $n\ge 1/4$, $|H_n(x)| \approx n^{-1/3}$ in $n \le x \le n + n^{1/3}$ and
$|H_n(x)| \approx x^{-1/4}(x-n)^{-1/4}$ in $x \ge n + n^{1/3}$.  These estimates
along with the formula for $\de_x\th(x,n)$ above yield \eqref{easy2}.

To prove Theorem \ref{nonclass},  we require symbol
estimates for $\th(x,n)$ expressed in terms of 
the variables $(\b,\n)$ with $x = \n \sec\b$, $n = \n$.

\begin{lem}\label{arg} Denote
\[
\s(\n,\b) = \th(\n\sec\b,\n)
\]
There is an absolute constant $C$ such that if 
$C\n^{-1/3} \le \b \le \pi/4$, then

a) $\disp |\de_\n^j\de_\b^k [\s(\b,\n) - \n(\tan\b - \b)] | \precsim
\n^{- j} \b^{-k} \quad j+k \ge 1$

b) $\disp |\de_\n^j\de_\b^k [\s(\b,\n)| \precsim
\n^{1- j} \b^{3-k} \quad j+k \ge 1 $
\end{lem}

In the transition region, where
$c\le \n\b^3 \le C$, for some absolute constants $0 < c < C < \infty$ 
more detailed asymptotics are required.  
Denote
$y = (\n^{2/3} \tan^2\b)/2$  and define the derivative of $\arg A(y)$ 
by 
\[
B(y) = \im \frac{A'(y)}{A(y)} = 
\im \frac{u(y) v'(y)  - u'(y) v(y)}{|A(y)|^2} 
\]
where $A(y) = u(y) + iv(y)$ is the Airy function in \eqref{wronskian}.
The asymptotic expansion for $\s(\n,\b)$ is given by
\begin{lem}\label{argasympt} If $j+k \ge 1$, then 

a) $\disp \de_\n \s = \frac13 (\n^{-1/3} \tan^2 \b) B(y)(1 + O(\b^2 + \n^{-4/3}\b^{-2}))$

b) $\disp \de_\b \s = (\n^{2/3} \tan \b  \sec^2 \b)  B(y)(1 + O(\b^2 +\n^{-1}\b^{-1}))$
\end{lem}

Proof of Lemma \ref{arg}.  First note that the definition of $\s$ and \eqref{adef} imply
\[
\de_\b \s = (\sec^2\b - 1)\n + \im \frac{\de_\b a}{a}
\]
and
\[
\de_\n \s = (\tan \b - \b) + \im \frac{\de_\n a}{a}
\]
Thus Lemma \ref{arg} follows from Proposition \ref{abound} and Corollary \ref{alowerbound}.

Next, we prove Lemma \ref{argasympt}.
\begin{align*}
\de_\n\s 
&= \tan \b -\b + \im (\de_\n a/a)  \\
&= \tan \b - \b + \im (\de_\n b/b)  
+ \im [(\de_\n a-\de_\n b)/a]  +  \im [(\de_\n b)(1/a - 1/b)]  \\
& = \tan\b - \b + \im (\de_\n b/b) + O(|\de_\n(a-b)|/|a| + |\de_n b||a-b|/|ab|) \\
& = \tan \b - \b + \im (\de_\n b/b)  + O((\b + \n^{-1/3})^5)
\end{align*}
Moreover (with $y =  (\n^{2/3}\tan^2\b)/2$), and recalling that $B(y) = \im (A'(y)/A(y))$, 
\[
\im \frac{\de_\n b}{b}  
= -\tan^3\b/3 + \frac{\n^{-1/3}}{3} \tan^2\b B(y)
\]
and $\tan \b - \b - (\tan^3\b)/3 \precsim \b^5$. Thus,
\[
\de_\n\s = \frac{\n^{-1/3}}{3} \tan^2\b 
B(y) + O((\b + \n^{-1/3})^5)
\]

Similarly, 
\begin{align*}
\de_\b\s 
&= \n(\sec^2\b - 1) + \im (\de_\b a/a)  \\
&= \n[\sec^2 \b - 1 + \im (\de_\b b/b) 
  + \im [(\de_\b a-\de_\b b)/a]
+ \im [(\de_\b b)(1/a - 1/b)] \\
& = \n[\sec^2 \b - 1]
 + \im (\de_\b b/b)   + O(\n (\b + \n^{-1/3})^4)
\end{align*}
Moreover,
\[
\im \frac{\de_\b b}{b}  
= -\n \tan^2\b\sec^2\b + \n^{2/3}(\tan\b)(\sec^2\b)
B(y)
\]
Since $\n[\sec^2 \b - 1] - \n \tan^2\b \sec^2\b \precsim \n \b^4$, 
\[
\de_\b\s   = 
\n^{2/3}(\tan\b)(\sec^2\b)B(y) + O(\n (\b + \n^{-1/3})^4)
\]
Finally,  in order to write the error as a multiplicative expression,
we use $|B(y)| \approx (1 + y)^{1/2}$ which follows from
the upper and lower bounds on $A(y)$ and $A'(y)$ (and
ultimately from the Wronskian formula).

We can now deduce Theorem \ref{nonclass} from Lemmas \ref{arg} and \ref{argasympt}.
Consider the functions $\n$ and $\b$ of $n$ and $m$ defined by
\[
\n = n;\quad \cos \b = n/\rho(m,n)
\]
For $3/4 \le m \le 3n$, 
$\r(m,n) -n \approx m^{2/3}n^{1/3}$ implies
\[
\b(m,n) \approx (m/n)^{1/3}\quad \mbox{if} \ 3/4 \le m \le 3n
\]
\[
\n \b^2 \approx m^{2/3}n^{1/3} \implies \b \approx (m/n)^{1/3}
\]
In particular when $m \ge 3/4$, $\b \ge c n^{-1/3}$ for some small absolute
constant $c>0$.   We will be treating three ranges of $\b$ in different
ways.  One is the transition region where 
$cn^{-1/3} \le \b \le Cn^{-1/3}$ where $C$ is a large constant.
The second is the region $Cn^{-1/3} \le \b \le \b_0$ for some small
absolute constant $\b_0$, and the third is the classical region
$\b_0 \le \b$.  We will never need to consider $\b$ smaller than $cn^{-1/3}$.

Differentiating $\s(\n,\b) = m\pi - \pi/2$, we find the implicit formulas
\[
\de_n \b = \frac{-\de_\n \s}{\de_\b \s};
\quad
\de_m \b = \frac{\pi}{\de_\b\s}
\]
We are going to prove by induction the property $Q(j,k)$ 
that says that for all $0 \le \ell_1 \le j$ 
and all $0 \le \ell_2 \le k$, 
\[
|\de_m^{\ell_1}\de_n^{\ell_2} \b| \precsim
\b^{1 - 3\ell_1}\r^{-\ell_1 - \ell_2} (\approx m^{1/3 - \ell_1}n^{-1/3 - \ell_2})
\]
($ \b \approx (m/n)^{1/3}$)
The property $Q(0,0)$ is trivial.   For $Q(0,1)$, observe that by Lemma \ref{argasympt}
\begin{equation}\label{ratio}
\frac{\de_\n \s}{\de_\b \s} = 
\frac{\n^{-1/3}\tan^2\b}{3\n^{2/3}\tan\b\sec\b}(1+ O(\b^2+ \n^{-4/3}\b^{-2}))
\end{equation}
Moreover, 
\[
\de_n\b = O(\b \n^{-1});\quad \de_m \b = O(\b^{-2}\n^{-1})
\]
The proof of $Q(1,0)$ is similar.  
Suppose that $Q(j,k)$ is valid.  Differentiating the implicit
formula for $\de_m^j \de_n^k \b$ with respect to $m$, there are
three types of things that can happen.  First the derivative
falls on the denominator, which is a power of $\de_\b \s$, in which case
the expression is multiplied by a constant times
\[
\frac{(\de_\b^2 \s)\de_m\b}{\de_\b s}
\]
But recall that $|\de_\b\s| \succsim \b^2\n$, $|\de_\b^2 \s| \precsim \b \n$,
and $|\de_m \b| \precsim \b^{-2}\n^{-1}$ so that the product is
majorized by $\b^{-3}\n^{-1}$.  This is the new factor required for
the estimate $Q(j+1,k)$.  If the derivative falls on the numerator,
then it may increase the degree of differentiation on a derivative of
$\b$, but always below the level of the induction hypothesis.  Replacing
a derivative of $\b$ by one derivative higher yields a change in estimation 
of the whole by the appropriate factor $\b^{-3}\n^{-1}$.  Finally,
the differentiation may land on a derivative of $D\s$.  This
replaces $D\s$ by $(\de_\b D\s) \de_m\b$, so the change in the 
estimation is the same (difference between $D\s$ and $\de_\b D\s$ is
a factor $\b^{-1}$ and $\de_m\b$ is bounded by $\b^{-2}\n^{-1}$.   Again
the product is $\b^{-3}\n^{-1}$, which is the factor we want.
Similarly, to prove $Q(j,k+1)$, differentiation with respect to $n$
produces an estimate that differs by the factor $\n^{-1}$ from the
bound for $Q(j,k)$. 

Next we can use $Q(j,k)$ to prove the property $P(j,k)$  that
for all $0 \le \ell_1 \le j$ and all $0 \le \ell_2 \le k$, 
\[
|\de_m^{\ell_1}\de_n^{\ell_2} (\r(m,n) - n)| \precsim
\b^{2 - 3\ell_1}\r^{1-\ell_1 - \ell_2} (\approx m^{2/3 - \ell_1}n^{1/3 - \ell_2})
\]
In fact, implicit differentiation of $\r \cos \b = n$ with respect to $n$ and with
respect to $m$ gives
\[
\de_n \r -1 = (\sec\b -1) + \r(\tan\b)\de_n \b
\]
and 
\[
\de_m \r = \r (\tan\b) \de_m\b
\]
The induction argument is similar to the proof of $Q(j,k)$ and is left to the 
reader.  This concludes Theorem \ref{nonclass} (a).  

We turn now to the lower bound (Theorem \ref{nonclass} b).  
\begin{align*}
\de_n \r -1 &= (\sec\b -1) + \r(\tan\b)\de_n \b \\
& = (\sec\b-1) - \r(\tan\b)\frac{\de_\n \s}{\de_\b \s} \\
& = \frac{1}2\b^2 + O(\b^4) - \frac13 \sin\b \tan\b (1 + O(\b^2 + \n^{-4/3}\b^{-2}) \\
& = \frac{1}6\b^2 + O(\b^4 + \n^{-4/3})
\end{align*}
Note that we need the precise asymptotics because the lower bound of order 
$\b^2 \approx m^{2/3}n^{-2/3}$ requires the coefficient $1/2 - 1/3  = 1/6 > 0$.  
The error term is lower order if $\b<<1$ and $n^{-4/3} << \b^2$.
Thus we have Theorem \ref{nonclass} (b) for $n$ sufficiently large,
in which case $ \n^{-2/3} <  < \b  \le C \n^{-1/3} $.  We do not need smaller
values of $\b$, because for $m\ge 1/20$, and $\b$ defined
implicitly in terms of $\r(m,n)$, $\b \approx (m/n)^{1/3} > > \n^{-2/3}$.

At last, here are some details of the much simpler estimates
in what we are calling the classical region.    
Let $n\ge 1$ and $x \ge (1+c)n$.  Then we are in the range in which
$(x-n) \sim x$.  As suggested in the remark of \cite{IJ} 
p. 1069, following the methods above one finds 
\[
H_n(x) = x^{-1/4}(x-n)^{-1/4}e^{i(x\sin\b - n\b)}s(x,n)
\]
where $\b(x,n)$ is defined by $x\cos \b = n$ and $s(x,n)$ satisfies
\[
|\de_x^j \de_n^k s(x,n)| \precsim x^{-j-k}
\]
Moreover, it is well-known that $|s(x,n)|>c>0$.  
It follows that $\th(x,n) = \arg H_n(x)$ (defined using any appropriate branch)
satisfies
\[
\th(x,n) = x\sin \b -n \b + E(x,n)
\]
where $E(x,n)$ is a symbol satisfying
\[
|\de_x^j \de_n^k E(x,n)| \precsim x^{-j-k}
\]
It is not hard to extend this estimate to negative values of $n$.  
Indeed the remaining range of $n$ is the range $-1/2 \le n \le 1$.
For a fixed range of the parameter $n$, the Hankel formula for $H_n(x)$ 
(\cite{W} pp. 196--198) may be used.
Let $x\ge 1$ and $-1/2 \le n \le x \cos \b_0$, where $\b_0$ is 
a small, fixed constant.   Then since $\b_0>0$,  
$x \ge (1+c)n$ for some $c>0$ and we are in the range
in which $(x-n)\simeq x$.  

One calculates that
\[
\de_x \b = \frac{n}{x^2 \sin \b}; \quad
\de_n\b = -\frac{1}{x\sin\b}
\]
and
\[
\de_n\th = -\b + \de_n E; \quad \de_x \th + \de_x E
\]
Differentiating the implicit equation for $\r(m,n)$
\[
\th(\r(m,n),n) = m\pi - \pi/2
\]
with respect to $m$ and $n$, we find
\[
\de_m\r = \frac{\pi}{\de_x\th(\r(m,n),n)} = \frac{\pi}{\sin\b} + F_1(\r,n)
\]
and
\[
\de_n\r = -\frac{\de_n\th(\r(m,n),n)}{\de_x\th(\r(m,n),n)} 
= \frac{\b}{\sin\b} + F_2(\r,n)
\]
with $F_1$ and $F_2$ satisfying 
\[
|\de_\r^j \de_n^k F(\r,n)| \precsim \r^{-1-j-k}
\]
and $\r\cos\b = n$, $\b_0 \le \b \le \pi/2 + \e$ ($\b$ extends a small amount 
beyond $\pi/2$ to accomodate the negative values of $n$ --- we only
care about $\r$ sufficiently large, so $\e$ can be arbitrarily small.)

Finally, this implicit asymptotic formula for $\r(m,n)$ can be differentiated
many times expressing derivatives of $\r(m,n)$ in terms of lower derivatives.
This yields 
\[
|\de_m^j \de_n^k \r(m,n)| \precsim \r^{1-j-k} \approx m^{1-j-k}
\]
for all $m \ge 100 + cn$ for a small fixed $c>0$.  
For the asymptotic gradient we already have the implicit formula
\[ 
(\de_m \r, \de_n\r) = (\pi/\sin\b, \b/\sin \b) + F(\r,n)
\]
with $F$ a symbol of order -1  as above.  Thus of order $O(1/m)$ 
since $\r\approx m$ in this range.  We also have the implicit
equation for $\r$ that gives us $\b$ as a function of $\r$ as follows:
\[
n(\tan \b - \b) + E(\r,n) = m\pi - \pi/2
\]
so that ($\r\cos\b = n$)
\[
\tan \b - \b = \frac{m\pi}{n}(1 - (E+\pi/2)/\pi m) = 
\frac{m\pi}{n}(1 - O(1/m))
\]
Thus if we define $\a(m,n)$ by
\[
\tan \a - \a =  \pi m/n
\]
(across $\pi/2$ as needed is perfectly ok) then
$\a = \b + O(1/m)$ and we get the asymptotic formula
we wanted.  (We can also characterize the error term as a symbol
rather than with a bound, but we don't need this.)

\section{Final Remarks}

The length spectrum of the disk is far from generic.  The regular
$k$-gon trajectory on the disk can be rotated around the circle giving
a one-parameter family of periodic geodesics of the same length,
whereas in general the length of the geodesic with $k$ reflections
varies depending on where its vertices are.  Nevertheless, we expect
that the contribution to the wave trace at times $t$ greater or equal 
to the perimeter from this family of $k$-reflection periodic trajectories
resembles the effect of the family of regular $k$-gons in the case of
the circle.  Indeed the almost integrable behavior of the
dyamical system of geodesic flow on convex planar regions near the
generalized geodesic that follows the boundary has been treated in
detail by Lazutkin \cite{L} and by Melrose and Marvisi \cite{MM}.  
Moreover, as we saw above, the asymptotics depended in a fundamental
way on nondegeneracy conditions on derivatives of Airy functions.  This 
is a hopeful sign, since general microlocal constructions of the
parametrix of the wave equation near the boundary also involve Airy
functions \cite{AM, MS, MT}.

As already mentioned, certain symbol properties of Bessel functions $J_\n(x)$ 
as a function of the two variables $x$ and $\n$ were already
proved in Ionescu and Jerison in \cite{IJ}.
There are three features of the treatment in \cite{IJ}
that are not enough for our purposes here. First, there is the minor
point that we need bounds on derivatives of all orders, not just the 
first two.  Second, the full symbol-type estimates (conjectured 
in \cite{IJ}, p. 1069) while valid 
are not sufficient here.  Those bounds are for
derivatives in the variables $(x,\n)$, whereas for the present purpose
in the range $\n\le x \le 2\n$, especially when $x$ is near $\n$,
it is necessary to distinguish a special directional derivative,
namely, the derivative with respect to $\n$ with the ratio
$x/\n$ held fixed (or equivalently with  $\b$ defined by $\cos \b =
\n/x$ held fixed).  Third, \cite{IJ} treats only upper 
bounds on $J_\n(x) + iY_\n(x)$ and its derivatives.  We need both upper and 
lower bounds on the argument and its derivatives, which require
detailed asymptotics, not just bounds on the Bessel functions and 
their first derivatives.


\begin{thebibliography}{9999}

\bibitem[AM]{AM} K. G. Andersson and R. B.  Melrose, The propagation of singularities
along gliding rays, \emph{Invent. Math.} \textbf{41} (1977) no. 3, 97--232.


\bibitem[C]{C} J. Chazarain, Construction de la parametrix du probl\`eme mixte
hyperbolique pour l'equation des ondes, \emph{C. R. Acad. 
Sci. Paris Ser. A-B} \textbf{276} (1973), A1213--1215.


\bibitem[CdV]{CdV} Y. Colin de Verdi\'ere, Spectre du laplacien et longueurs
des geodesiques periodiques II, Compositio Math. \textbf{27} (1973), 159--184.


\bibitem[DG]{DG} J. J. Duistermaat and V. Guillemin, 
The spectrum of positive elliptic operators and periodic bicharacteristics,
\emph{Invent. Math.} \textbf{29} (1975) 39--79.

\bibitem[GM]{GM} V. Guillemin and R. B.  Melrose, The Poisson Summation formula
for manifolds with boundary, \emph{Adv. in Math.} \textbf{32} (1979) 204--232.

\bibitem[HeZ]{HeZ} H. Hezari and S. Zelditch, Inverse spectralproblem
for analytic $(\Z/2\Z)^n$-symmetric domains in $\R^n$, arXiv: 0902.1373 

\bibitem[IJ]{IJ} 
A. D. Ionescu and D. Jerison,  
On the absence of positive eigenvalues of Schr\"odinger operators with 
rough potentials,
\emph{Geom. Funct. Anal.}  \textbf{13}  (2003),  no. 5, 1029–1081. 



\bibitem[L]{L} V. F.  Lazutkin, Construction of an asymptotic series of 
eigenfunctions of ``bouncing ball'' type, \emph{Proc. Steklov Inst. Math.}
(1968) 125--140.


\bibitem[MM]{MM} S. Marvisi and R. B. Melrose, 
Spectral invariants of convex planar regions, \emph{J. Diff. Geom.} \textbf{17}
(1982) 475--502.

\bibitem[MS]{MS} R. B. Melrose and J. Sj\"ostrand,
Singularties of boundary value problems,
\emph{Comm. Pure Appl. Math.}
\textbf{31}
(1978) 593--617.


\bibitem[MT]{MT} R. B. Melrose and M. Taylor,
Near peak scattering and the corrected Kirchkoff approximation for a  convex
obstacle,
\emph{Adv. in Math.} 
\textbf{55}
(1985) 242--315.



\bibitem[W]{W} G. N. Watson,
\emph{A treatise on the theory of Bessel functions}, 2nd edition,
Cambrige U. Press  (1966)

\bibitem[Z]{Z} S. Zelditch, 
Inverse spectral problems for analytic domains II: domains with one symmetry,
\emph{Ann. of Math.} (2) 
\textbf{170}
(2009)  205--269.



\end{thebibliography}
\end{document}